\documentclass{jgcc} 

 \pdfoutput=1
\usepackage{lastpage}
\jgccdoi{18}{1}{4}{16639}
\jgccheading{}{\pageref{LastPage}}{}{}{Oct.~5,~2025}{Jun.~17,~2026}{}

\keywords{free group, automorphisms, automorphic orbits}

\usepackage{hyperref}
\usepackage{orcidlink}

\begin{document}

\title{Automorphic orbits in free groups: recent progress}

\author[V.~Shpilrain]{Vladimir Shpilrain\, 
 \orcidlink{0000-0003-2239-7996} 
}	
\address{Department of Mathematics, The City College of New York, New York, NY 10031, USA}	
\email{shpilrain@yahoo.com}  
\NewDocumentCommand{\Aut}{}{\operatorname{Aut}}
\NewDocumentCommand{\Wh}{}{\operatorname{Wh}}
\NewDocumentCommand{\defnsub}{}{\textit}

\NewDocumentCommand{\CCLR}{}{strongly reduced}
\NewDocumentCommand{\RCR}{}{weakly reduced}
\NewDocumentCommand{\Rone}{}{partially reduced}

\newcommand{\Z}{\mathbb{Z}}

\newcommand{\Q}{\mathbb{Q}}

\newcommand{\F}{\mathbb{F}}

\newcommand{\N}{\mathbb{N}}

\newcommand{\R}{\mathbb{R}}

\newcommand{\C}{\mathbb{C}}

\begin{abstract}
In this survey, we describe recent progress on asymptotic properties of various automorphic orbits in free groups. In particular, we address the problem of counting potentially positive elements of a given length. We also discuss complexity (worst-case, average-case, and generic-case) of Whitehead's automorphism problem and relevant properties of automorphic orbits, including orbit-blocking words.
\end{abstract}

\maketitle

\hfill {\small To Alexei Myasnikov on his special birthday}

\section{Introduction}\label{intro}

In this survey article, $F_r$ will denote a free group of rank $r$.
Free groups and their properties are at the core of combinatorial group theory; they have been studied for over 100 years now.
Automorphisms of free groups, in particular, proved to be attractive objects of research, bridging group theory, topology, geometry, and dynamics.

Given an element $u \in F_r$, by the {\it automorphic orbit} of $u$ we mean the set of all automorphic images of $u$, i.e., the set $\{\varphi(u)\}$ where $\varphi$ runs over all automorphisms of $F_r$.
Automorphic orbits in free groups were studied by J.~H.~C.~Whitehead in his now classical papers \cite{Wh} and \cite{Wh2}.

Since then, there was a vast amount of research on various properties of automorphic orbits in $F_r$. We can single out the research avenue focused on fixed points of automorphisms that was popular in the late 20th century and culminated in a seminal paper by Bestvina and Handel \cite{BHfix}.

The early 21st century saw a surge in interactions between group theory and theoretical computer science. Most of these interactions are well covered in the monograph \cite{bookof10}, but one particular trend that stands out is emerging interest in complexity of group-theoretic algorithms. This includes not only the ``traditional" worst-case complexity but also generic-case and average-case complexity introduced into group theory in \cite{KMSS1} and \cite{KMSS2}, respectively.

In this article, we briefly survey research directions that stemmed from the paper \cite{MS}. The main theme is complexity of Whitehead's problem (Section \ref{WP}) that led to exciting progress on orbit-blocking words (Section \ref{blocking}). Then, there are automorphic orbits of potentially positive elements (Section \ref{pp}); these, too, saw remarkable progress recently.

In Section \ref{Trans}, we consider pairs of elements of $F_r$ that have ``similar" automorphic orbits; elements like that are called (boundedly) translation equivalent following \cite{KLSS}. Somewhat surprisingly, there was no progress on characterizing or recognizing (boundedly) translation equivalent elements in the last 15 years or so, but we include this problem here in the hope that the right tools to attack it will eventually be found.

Another interesting direction (not covered in the present survey) is studying intersection of the automorphic orbit of a given element of $F_r$ with a given subgroup or a rational subset of $F_r$, see \cite{Silva}. This is also relevant to potential positivity, see our Section \ref{pp}.

Finally, we note that sometimes the words ``automorphic orbit" mean something different (smaller), see e.g. \cite{Brinkmann} or \cite{Ventura}. Namely, given an element $u \in F_r$ and an automorphism $\varphi$ of $F_r$, an automorphic orbit is defined as the set $\{\varphi^n(u), n=1, 2, \ldots\}$. This, too, leads to various interesting questions that are not in the scope of the present survey.

\section{Orbit-blocking words}\label{blocking}

Let $F_r$ be a free group generated by $x_1,\ldots, x_r$.
For a given element $w \in F_r$, we say that a word $v$ is $w$-orbit-blocking (or just orbit-blocking if $w$ is clear from the context) if for any automorphism $\varphi \in \Aut(F_r)$, the word $v$ is not a subword of the cyclic reduction of $\varphi(w)$.

{\it Primitivity-blocking} words are the same as $x_1$-orbit-blocking.
The existence of primitivity-blocking words easily follows from Whitehead's observation that the Whitehead graph of any cyclically reduced primitive element is either disconnected or has a cut vertex.

The {\it Whitehead graph} $\Wh(w)$ of $w \in F_r$ has $2r$  vertices that correspond to the generators $x_i$ and their inverses.
For each occurrence of a subword $x_i x_j$ in the word $w \in F_r$, there is an edge in $\Wh(w)$ that connects the vertex $x_i$ to the vertex $x_j^{-1}$; if $w$ has a subword $x_i x_j^{-1}$, then there is an edge connecting $x_i$ to $x_j$, etc.
There is one more edge (the external edge): this is the edge that connects the vertex corresponding to the last letter of $w$ to the vertex corresponding to the inverse of the first letter.

It was observed by Whitehead himself in his \textit{cut vertex lemma} (see \cite{Wh2}) that the Whitehead graph of any cyclically reduced primitive element $w$ is either disconnected or has a cut vertex, i.e., a vertex that, having been removed from the graph together with all incident edges, increases the number of connected components of the graph.
A short and elementary proof of this result was recently given in \cite{Heusener}, and a more general case of primitive elements in subgroups of $F_r$ was recently treated in \cite{Ascari}.

Thus, for example, if the Whitehead graph of $w$ (without the external edge) is complete (i.e., any two vertices are connected by at least one edge), then $w$ is primitivity-blocking because in this case, if $w$ is a subword of $u$, then the Whitehead graph of $u$, too, is complete and therefore is connected and does not have a cut vertex.
Here are some examples of primitivity-blocking words: $x_1^n x_2^n \cdots x_r^n x_1$ (for any $n \ge 2$), $[x_1, x_2][x_3, x_4]\cdots [x_{n-1}, x_{n}]x_1^{-1}$ (for an even $n$), etc.
Here $[x, y]$ denotes $x^{-1}y^{-1}xy$.

In \cite{rank2}, it was shown that for any $w \in F_2$, there are $w$-orbit-blocking words.
This was based on an explicit description of bases (equivalently, of automorphisms) of $F_2$ from \cite{Cohen}. Then, in \cite{rankr}, the following general result was established:

\begin{thm}\label{Tblocking}
  Let $w\in F_r$ have cyclically reduced length $\ell$.
  Let $\{v_i\}_{i=1}^{\ell+1}$ be a sequence of primitivity-blocking words such that there is no cancellation between adjacent terms of the sequence, nor between $v_{\ell+1}$ and $v_1$.
  Then  $\prod_{i=1}^{\ell+1}v_i$ is $w$-orbit-blocking.
\end{thm}

Note that, interestingly, $w$-orbit-blocking words described in Theorem \ref{Tblocking} depend on the length of $w$ but not on the actual $w$.

As can be seen from the statement of Theorem \ref{Tblocking}, primitivity-blocking words play an important role, so we are going to discuss them separately.

\subsection{Primitivity-blocking words}

By definition, primitivity-blocking words are those that cannot be subwords of any {\it cyclically reduced} primitive elements.
What is needed to prove Theorem \ref{Tblocking} is words with somewhat stronger property: they cannot be subwords even of some primitive elements that are not cyclically reduced.

\begin{defi}
Call a basis of $F_r$ \defnsub{\CCLR} if there is no element $g \in F_r$ such that conjugating each element of the basis by $g$ decreases the sum of the lengths of the basis elements.
\end{defi}

Then there is the following technical result of independent interest.

\begin{lem}\cite{rankr}
  \label{lem:big boy}
  Let $B$ be a \CCLR\ basis of $F_r$.
 Then no primitivity-blocking word appears as a subword of any element in $B$.
\end{lem}

Even though Lemma \ref{lem:big boy} suffices to prove Theorem \ref{Tblocking}, it can be generalized further.

\begin{defi}
  Call a basis $B$ \defnsub{\RCR} if there is no letter $p$ such that every element of $B$ begins with $p$ and ends with $p^{-1}$.
\end{defi}

Then we have:

\begin{lem}\cite{rankr}
  \label{thm:RCR no pbsw}
  Let $B$ be a \RCR\ basis of $F_r$, where $r\geq 3$, and let $w \in B$.
Then $w$ does not contain any primitivity-blocking subwords.
\end{lem}

Note that in $F_2$, a \RCR\ basis can contain a primitivity-blocking word.
For instance, the basis $\{ x_1x_2x_1^{-1}, x_2x_1^{-1} \}$ is \RCR, but contains $x_1x_2x_1^{-1}$, which is primitivity-blocking in $F_2$.

Then, to show that blocking primitivity does not have to be due to connectivity properties of the Whitehead graph, we include the following examples.

\begin{prop}\cite{rankr}\label{exam} \noindent {\bf (a)}
In the group $F_2$, the words $x_1^k x_2^k$ are primitivity-blocking for any $k \ge 2$.

\noindent {\bf (b)} In the group $F_r$, $r \ge 3$, the words $x_1^k \ldots  x_r^k$ are not primitivity-blocking for any $k \ge 1$.
\end{prop}

Part (b) of Proposition \ref{exam} can be generalized as follows:

\begin{prop}\cite{rankr}\label{XY}
Let $F_r$ where $r \ge 3$ be generated by the set $X\sqcup Y$.
Let $\langle X \rangle$ and $\langle Y \rangle$ denote the subgroups of $F_r$ generated by $X$ and $Y$, respectively.
Let $w=w_Xw_Y$, where $w_X \in \langle X \rangle$, $w_Y \in \langle Y \rangle$.
Then $w$ is not a primitivity-blocking word.
\end{prop}

We note that the shortest  primitivity-blocking word in $F_2$ is $x_1^{-1}x_2x_1$.
The Whitehead graph (with or without the external edge) of this word does have a cut vertex though.

It is natural then to look for the shortest primitivity-blocking words in $F_r$ for $r>2$.

\begin{thm}\cite{rankr}\label{thm:shortpb}
  The word $w = x_1x_2x_3\dots x_{r-1}x_r^2x_{r-1}\dots x_3x_2x_1^{-1}\in F_r$ is the shortest primitivity-blocking word for $r>2$.
\end{thm}

A similar word $w = x_1x_2 \dots x_{r-1}x_r^2x_{r-1}\dots x_2x_1$ is not primitivity-blocking in any $F_r$, $r\ge 2$. Indeed, for $r=2$, $w$ is a subword of $x_2x_1x_2^2x_1$, which is a cyclically reduced primitive word.

For $r > 2$, let $u=x_2x_3\dots x_{r-1}x_r^2x_{r-1}\dots x_3x_2$.
Consider the word $v=ux_1ux_1x_3$. This  $v$ clearly contains $w$ as a subword and is cyclically reduced. Apply the automorphism $\varphi = x_1\mapsto u^{-1}x_1, x_i \mapsto x_i, ~i \ge 2$.
  Then $\varphi(v) = x_1^2x_3$ is a primitive element, and therefore so is $v$.

Finally, we mention that primitivity-blocking words in $F_2$ can be algorithmically recognized. For $F_r, r>2$, the problem is open.

\begin{thm}\cite{rankr}\label{F2alg}
    There is an algorithm that decides, for a given word $u \in F_2$ of length $n$, whether or not $u$ is primitivity-blocking and has the worst-case time complexity $O(n^2)$.
\end{thm}

\section{Whitehead's automorphism problem}\label{WP}

Let $F_r$ be a free group with a free generating set $x_1, \ldots, x_r$.
Whitehead's automorphism problem is: given $u, v \in F_r$, decide whether or not $v$ is an automorphic image of $u$. The problem is solved by {\it Whitehead's algorithm} that dates back to 1936 \cite{Wh}:
\medskip

(1) First reduce $u$ to a word of minimal possible length by applying  ``elementary" Whitehead automorphisms. The same for $v$. This part is ``greedy" and has quadratic time complexity with respect to the length of the input.
\medskip

(2) If the reduced words are not of the same length, then $u$ and $v$ are not in the same automorphic orbit. If they are of the same length, then things get interesting (in terms of complexity). Indeed, the procedure outlined in the original paper by Whitehead suggested this part of the algorithm to be of exponential time with respect to the length of
the words since it involves going through all images of a given $u$, having the same length as $u$ does, under elementary  Whitehead automorphisms. This raises a natural question about the maximum possible number of elements of length $m$ in the automorphic orbit of $u$.
\medskip

Thus, a challenging open problem is:

\begin{prob}\label{worst-case}
What is the worst-case time complexity of the second part of Whitehead's algorithm?
\end{prob}

In $F_2$, the worst-case time complexity was shown to be at most quadratic \cite{MS}.
Actually, if the length of $u \in F_2$ is $m$ and it cannot be reduced by any Whitehead automorphism, then the maximum possible number of elements of length $m$ in the automorphic orbit of $u$ is precisely $8m-40$ for $m \ge 10$, see \cite{Khan}.

In $F_3$, the maximum number is $48m^3 - 480m^2 + 1104m - 672$ for $m \ge 11$ as suggested (but not proved) by C. Sims based on computer experiments. A particular Whitehead-reduced word of length $m$ whose automorphic orbit (limited to elements of length $m$) has the cardinality given by the latter polynomial is
$u=x_1^k x_2 x_1x_2^{-1} x_1 x_2^2x_3^2,$ where $k=m-8$.

Let $A(u)$  denote the set of elements  $\{v \in F_r;  ~|v| = |u|,  ~\varphi(v)=u$  for some $\varphi \in Aut(F_r)\}$.
D. Lee \cite{Leeorbit} came close to solving Problem \ref{worst-case} completely. Specifically, she proved that the cardinality of $A(u)$  is bounded by a polynomial function of $|u|$ under the following condition: If two letters $x_i$ (or $x_i^{-1}$) and $x_j$ (or $x_j^{-1}$) with $i<j$ occur in $u$, then the total number of $x_i^{\pm 1}$ occurring in $u$ is strictly less than the total number of $x_j^{\pm 1}$ occurring in $u$.

Then, D. Lee \cite{Leetight} proved that, under the same assumption on $u$, the cardinality of $A(u)$ is bounded by a polynomial function of $m=|u|$ of degree $2r-3$, and that this bound is sharp.

We can therefore make Problem \ref{worst-case} somewhat more specific:

\begin{prob}
Let $u$ be a Whitehead-reduced word of length $m$, and let $A(u)$ be the set of automorphic images of $u$ that have length $m$. Is the cardinality of $A(u)$ bounded by a polynomial function of $m$?
\end{prob}

\subsection{Generic-case and average-case complexity of Whitehead's problem}\label{average}

Informally, generic-case complexity of an algorithm is complexity on ``most" inputs, or on ``random" inputs. For a formal treatment, see \cite{KMSS1}.

Generic-case complexity of the classical Whitehead's algorithm (see above) is linear, as shown in \cite{kssgen}. The reason is that for a word $u$ selected uniformly at random from the set of elements of length $n$ in $F_r$, applying any Whitehead's automorphism (except permutations on the set of generators and their inverses) increases the length of $u$ with probability $1-O(2^{-n})$.

To get a meaningful (e.g. subexponential) estimate of the average-case complexity of Whitehead's problem, one has to get a more precise (than just ``exponential") upper bound on the worst-case complexity. It may be then  possible to get a subexponential bound on the average-case complexity of Whitehead's problem without getting a subexponential bound on the worst-case complexity.

In \cite{Shpilrain2023}, \cite{rank2}, \cite{rankr}, the following variant of Whitehead's problem was considered: given a fixed $u \in F_r$, decide, on an input $v \in F_r$ of length $n$, whether or not $v$ is an automorphic image of $u$. Thus, in this variant the input consists of just one word.

It turns out that the average-case complexity of this version of the Whitehead problem is constant if the input $v$ is a cyclically reduced word.
For a formal definition of the average-case complexity of an algorithm in the context of group theory we refer the reader to \cite{KMSS2}.

The algorithm with constant average-case complexity from \cite{rankr} is a combination of two different algorithms running in parallel: one is fast but may be inconclusive, whereas the other one is conclusive but relatively slow.
This idea has been used for group-theoretic algorithms since at least \cite{KMSS2}.

Before running the two algorithms, a pre-computation would be performed on $u$.
Specifically, $u$ is reduced to an element of minimum length in its automorphic orbit. This takes worst-case quadratic time in $|u|$.
However, since $u$ is a fixed word, this amounts to constant time for the algorithm.
Denote the obtained element of minimum length in the orbit of $u$ by $\bar u$.
Once the word's length has been minimized, the two algorithms will be run in parallel on the result.

A fast algorithm $\mathcal{T}$ would detect a $\bar u$-orbit-blocking subword $B(\bar u)$ of a (cyclically reduced) input word $v$, as follows.
Let $n$ be the length of $v$.
The algorithm $\mathcal{T}$ would read the initial segments of $v$ of length $k$, $k=1, 2, \ldots,$ adding one letter at a time, and check if this initial segment has $B(\bar u)$ as a subword.
This takes time bounded by $C\cdot k$ for some constant $C$, see e.g.\ \cite[p.\ 338]{Knuth2}.

The ``usual" Whitehead algorithm, call it $\mathcal{W}$, would minimize $|v|$ taking time quadratic in $|v|$.
Denote the obtained element of minimum length in the orbit of $v$ by $\bar v$.
If $|\bar v| \ne |\bar u|$, then $\mathcal{W}$ stops and reports that $v$ is not in the automorphic orbit of $u$.
If $|\bar v| = |\bar u|$, then the algorithm $\mathcal{W}$ would apply all possible sequences of elementary Whitehead automorphisms that do not change the length of $\bar v$ to see if any of the resulting elements are equal to $\bar u$.
This part may take exponential time in $|\bar v| = |\bar u|$, but since we consider $|u|$ constant with respect to $|v|$ and $|u|$ bounds $|\bar u|$ above, exponential time in $|\bar u|$ is still constant with respect to $|v|$.
Thus, the total time that the algorithm $\mathcal{W}$ takes is quadratic in $|v|$.

Finally, the algorithm $\mathcal{A}$ will consist of the algorithms $\mathcal{T}$ and $\mathcal{W}$ running in parallel.
Then we have:

\begin{thm}\label{average-case}\cite{rankr}
Suppose possible inputs of the above algorithm $\mathcal{A}$ are cyclically reduced words that are selected uniformly at random from the set of cyclically reduced words of length $n$.
Then the average-case time complexity (i.e., the expected runtime) of the algorithm $\mathcal{A}$, working on a classical Turing machine, is $O(1)$, a constant that does not depend on $n$.
If one uses the ``Deque" (double-ended queue) model of computing \cite{deque} instead of a classical Turing machine, then the ``cyclically reduced" condition on the input can be dropped.
\end{thm}

\section{Potentially positive elements}\label{pp}

An element $u$ of $F_r$ is called {\it positive} if no $x_i$ occurs in   $u$   to a negative exponent. An element  $u$ is called {\it  potentially positive}  if $\alpha(u)$ is positive for some automorphism $\alpha$ of the group $F_r$.

An element  $u$ is called {\it stably potentially positive}  if $\alpha(u)$ is positive for some automorphism $\alpha$ of the group $F_n$ for some $n \ge r$. (Here we tacitly assume that $F_r$ is canonically embedded in $F_n$.)

Motivation for considering potentially positive elements comes from the fact that various properties of one-relator groups are easier to establish if the relator is a positive word. For example, Baumslag \cite{Baumslag} showed that one-relator groups with a positive relator are residually solvable, and Wise \cite{Wise} showed that one-relator groups with a positive relator (satisfying a small cancellation condition) are residually finite. All these properties obviously hold upon replacing ``positive" with ``potentially positive" or even
``stably potentially positive".

Clark and Goldstein \cite{Clark} proved that all stably potentially positive elements are potentially positive.

Goldstein \cite{Goldstein} and D. Lee \cite{Leepositive} offered (exponential time) algorithms for deciding potential positivity in $F_2$. Another algorithm was offered by Silva and Weil \cite{Silva}.

Recently, Koch-Hyde, O'Connor, and Olive \cite{pprank2} reported an algorithm with the worst-case time complexity $O(n^2)$. This algorithm has linear generic-case complexity. If inputs are cyclically reduced, then the average-case complexity of this algorithm is constant. This is because, as shown in \cite{pprank2},  there are potential positivity-blocking words in $F_2$, i.e., words that cannot be subwords of any potentially positive word. An example would be $xyx^{-1}y^{-1}x$.

A very interesting and challenging question is: how many potentially positive elements of length $n$ are there in $F_r$?

The set of positive elements is exponentially negligible in the set of potentially positive elements because any element with only positive occurrences of all but one of the generators is potentially positive. (Note that the number of positive elements of length $n$ in $F_r$ is $r^n$, whereas the total number of elements of length $n$ is $2r\cdot (2r-1)^n$.)

The number of potentially positive elements of length $n$ is bounded from below by the number of elements of length $n$ with only positive occurrences of all but one of the generators. The latter number is $> (\frac{r^2+2r-1}{r+1})^n$.

For $r>3$, a larger lower bound is provided by the number of primitive elements of length $n$, which is $O(n (2r-3)^n)$ by \cite{PuderWu}.

Recently, \cite{ppositive} obtained a tight estimate for the growth rate in the case $r=2$ by proving that the number of potentially positive elements of length $n$ in $F_2$ is $O((\lambda + \epsilon)^n)$ for any $\epsilon > 0$, where $\lambda \approx 2.505$ is the largest root of the polynomial $\lambda^4 - 3 \lambda^3 +\lambda^2 + \lambda -1$.

Summing up, we ask:

\begin{prob}
Let $P(r, n)$ denote the number of potentially positive elements of length $n$ in $F_r$. What is the growth of $P(r, n)$ as a function of $n$, with $r$ fixed ?
\end{prob}

\section{Translation equivalence}\label{Trans}

Two elements  $u, v$ of a free group $F_r$ are called {\it translation equivalent} if for every free and discrete isometric action of $F_r$ on an $\R$-tree, translation lengths of $u$ and $v$ are equal.

A purely combinatorial characterization of translation equivalence was given in \cite{KLSS}: elements  $u, v$ are translation equivalent if the cyclic length of $\varphi(u)$ equals the cyclic length of $\varphi(v)$ for every automorphism $\varphi$ of the group $F_r$.

A ramification of this definition is: two elements  $u, v$  are called {\it boundedly translation equivalent} if the ratio of the cyclic lengths of $\varphi(u)$ and  $\varphi(v)$ is bounded away from 0 and from $\infty$ when $\varphi$ runs through all automorphisms of $F_r$.

It is natural to ask:

\begin{prob}
Is there an algorithm which, when given two elements of a finitely generated free group, decides whether or not they are (boundedly) translation equivalent?
\end{prob}

It is not quite trivial to even produce examples of (boundedly) translation equivalent elements. In \cite{KLSS}, it was shown that in $F_2$, any $u(x,y)$ is translation equivalent to $u(x,y)$ read backwards.

Another source of translation equivalent elements comes from {\it character equivalence}. We say that $u$ and $v$ are {\it trace equivalent} or {\it character equivalent} if for every representation $\alpha: F_r \to SL_2(\C)$ one has $tr(\alpha(u))=tr(\alpha(v))$. Character equivalent words come from the so-called ``trace identities” in $SL_2(\C)$ and are quite plentiful (see, for example, \cite{Horowitz}).

D. Lee showed \cite{LeeTranslation} that whenever $g$ is translation equivalent to $h$ in $F_r$ and $w(x,y)$ is arbitrary, one has $w(g,h)$ translation equivalent to $w(h,g)$ in $F_r$. She also reported an algorithm that decides translation equivalence in $F_2$.

In \cite{LeeTranslationrank2}, D. Lee offered  an algorithm that decides bounded translation equivalence in $F_2$.

Surprisingly, there has not been any progress on (bounded) translation equivalence since then.

\end{document}